\documentclass[journal,12pt,onecolumn,draftclsnofoot]{IEEEtran}
\usepackage[latin9]{luainputenc}
\usepackage{float}
\usepackage{calc}
\usepackage{amsmath}
\usepackage{amsthm}
\usepackage{graphicx}
\usepackage[unicode=true]
 {hyperref}

\makeatletter
 \let\oldforeign@language\foreign@language
 \DeclareRobustCommand{\foreign@language}[1]{%
   \lowercase{\oldforeign@language{#1}}}
\theoremstyle{plain}
\newtheorem{thm}{\protect\theoremname}
\theoremstyle{plain}
\newtheorem{lem}[thm]{\protect\lemmaname}

\usepackage[caption=false,font=footnotesize]{subfig}
\usepackage{algorithmic}
\usepackage{cite}
\usepackage{amssymb}

\usepackage{graphicx}
\usepackage{epsfig}
\usepackage{epstopdf}

\makeatother

\providecommand{\lemmaname}{Lemma}
\providecommand{\theoremname}{Theorem}

\begin{document}

\title{MIMO Transmit Beampattern Matching Under Waveform Constraints}

\author{Ziping~Zhao,~\IEEEmembership{Student Member,~IEEE}, and Daniel~P.~Palomar,~\IEEEmembership{Fellow,~IEEE}\thanks{This work was supported by the Hong Kong RGC 16206315 research grant.
The author Z. Zhao is supported by the Hong Kong PhD Fellowship Scheme
(HKPFS).}\thanks{The authors are with the Department of Electronic and Computer Engineering,
The Hong Kong University of Science and Technology (HKUST), Clear
Water Bay, Kowloon, Hong Kong (e-mail: \protect\href{mailto:ziping.zhao@connect.ust.hk}{ziping.zhao@connect.ust.hk};
\protect\href{mailto:palomar@ust.hk}{palomar@ust.hk}).}}

\markboth{\MakeLowercase{}}{ZHAO AND PALOMAR: MRP with Budget Constraint}
\maketitle
\begin{abstract}
In this paper, the multiple-input multiple-output (MIMO) transmit
beampattern matching problem is considered. The problem is formulated
to approximate a desired transmit beampattern (i.e., an energy distribution
in space and frequency) and to minimize the cross-correlation of signals
reflected back to the array by considering different practical waveform
constraints at the same time. Due to the nonconvexity of the objective
function and the waveform constraints, the optimization problem is
highly nonconvex. An efficient one-step  method is proposed to solve
this problem based on the majorization-minimization (MM) method. The
performance of the proposed algorithms compared to the state-of-art
algorithms is shown through numerical simulations.
\end{abstract}


\section{Introduction}

Multiple-input multiple-output (MIMO) systems \cite{LiStoica2008}
have the capacity to transmit independent probing signal or waveforms
from each transmit antenna. Such waveform diversity feature leads
to many desirable properties for MIMO systems. For example, a modern
MIMO radar has many appealing features, like higher spatial resolution,
superior moving target detection and better parameter identifiability,
compared to the classical phased-array radar \cite{BlissForsythe2003,FishlerHaimovichBlumChizhikCiminiValenzuela2004,ForsytheBlissFawcett2004}. 

The MIMO transmit beampattern matching problem is critically important
in many fields, like in defense systems, communication systems, and
biomedical applications. This problem is concerned with designing
the probing waveforms to approximate a desired antenna array transmit
beampattern (i.e., an energy distribution in space and frequency)
and also to minimize the the cross-correlation of the signals reflected
back from various targets of interest by considering some practical
waveform constraints. The MIMO transmit beampattern matching problem
appears to be difficult from an optimization point of view because
the existence of the fourth-order nonconvex objective function and
the possibly nonconvex waveform constraints which are used to represent
desirable properties and/or enforced from an hardware implementation
perspective \cite{Skolnik1990}. 

In \cite{FuhrmannSanAntonio2004}, the MIMO transmit beampattern matching
problem was formulated to minimize the difference between the designed
beampattern and the desired one. The formulation in \cite{FuhrmannSanAntonio2004}
was modified in \cite{FuhrmannSanAntonio2008,StoicaLiXie2007} by
introducing the cross-correlation between the signals. And in \cite{StoicaLiXie2007},
the authors proposed to design the waveform covariance matrix to match
the desired beampattern through semidefinite programming. A closed-form
waveform covariance matrix design method was also proposed based on
discrete Fourier transform (DFT) coefficients and Toeplitz matrices
in \cite{LiporAhmedAlouini2014,BouchouchaAhmedAl-NaffouriAlouini2017}.
But such kind of methods can perform badly for small number of antennas.
After the waveform covariance matrix is obtained, other methods should
be applied to synthesize a desired waveform from its covariance matrix.
For example, a cyclic algorithm was proposed in \cite{StoicaLiZhu2008}
to synthesize a constant modulus waveform from its covariance matrix.
These methods are usually called two-steps methods. In practice, they
could become inefficient and suboptimal if more waveform constraints
are considered.

In \cite{WangWangLiuLuo2012}, it was found that directly designing
the waveform to match the desired beampattern can give a better performance,
which is referred to as the one-step method. But the method in \cite{WangWangLiuLuo2012}
is tailored to the constant modulus constraint and can be slow in
convergence. In \cite{ChengHeZhangLi2017}, the problem was solved
based on the alternating direction method of multipliers (ADMM) \cite{BoydParikhChuPeleatoEckstein2011}.
However, again the proposed algorithm is only designed for dealing
with unimodulus constraint.

The majorization-minimization (MM) method \cite{HunterLange2004,SunBabuPalomar2016}
has shown its great efficiency in deriving fast and convergent algorithms
to solve nonconvex problems in many different applications \cite{SongBabuPalomar2015,ZhaoPalomar2018}.
In this paper, we propose a one-step method to directly solve the
MIMO transmit beampattern matching problem based on the MM method
by considering different waveform constraints. The performance of
our algorithms compared to the existing algorithms is verified through
numerical simulations.

\section{MIMO Transmit Beampattern Matching Problem Formulation}

\begin{figure}[t]
\centering{}\includegraphics[scale=0.7]{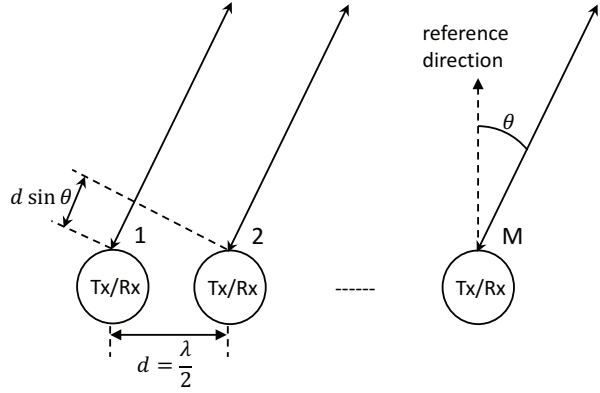}\caption{\label{fig:MIMO-radar-transceiver}MIMO transceiver with $M$ antennas
and $\theta$ is the spacial direction of interest.}
\end{figure}

A colocated MIMO radar \cite{LiStoica2007} with $M$ transmit antennas
in a uniform linear array (ULA), as shown in Fig. \ref{fig:MIMO-radar-transceiver},
is considered. Each transmit antenna can emit a different waveform
$x_{m}\left(n\right)$ with $m=1,2,\ldots,M$, $n=1,2,\ldots,N$,
where $N$ is the number of samples. Let $\mathbf{x}\left(n\right)=\Bigl[x_{1}\left(n\right),x_{2}\left(n\right),\ldots,x_{M}\left(n\right)\Bigr]^{T}$
be the $n$th sample of the $M$ transmit waveforms and $\mathbf{x}=\Bigl[\mathbf{x}^{T}\left(1\right),$
$\mathbf{x}^{T}\left(2\right),\ldots,\mathbf{x}^{T}\left(N\right)\Bigr]^{T}$
denote the waveform vector.

The signal at a target location with angle $\theta$ ($\theta\in\Theta$,
which is the angle set) is represented by
\[
\sum_{m=1}^{M}e^{-j\pi\left(m-1\right)\sin\theta}x_{m}\left(n\right)=\mathbf{a}^{T}\left(\theta\right)\mathbf{x}\left(n\right),\,n=1,\ldots,N,
\]
where $\mathbf{a}\left(\theta\right)$ is the transmit steering vector
written as $\mathbf{a}\left(\theta\right)=\left[1,e^{-j\pi\sin\theta},\ldots,e^{-j\pi\left(M-1\right)\sin\theta}\right]^{T}$.
Then, the power for the probing signal $\mathbf{x}$ at location $\theta$
which is named the \emph{transmit beampattern} can be written as follows:
\[
\begin{aligned} & P\left(\theta,\mathbf{x}\right)\\
= & \sum_{n=1}^{N}\left(\mathbf{a}^{T}\left(\theta\right)\mathbf{x}\left(n\right)\right)^{\ast}\left(\mathbf{a}^{T}\left(\theta\right)\mathbf{x}\left(n\right)\right)\\
= & \left(\left(\mathbf{I}_{N}\otimes\mathbf{a}^{T}\left(\theta\right)\right)\mathbf{x}\right)^{H}\left(\left(\mathbf{I}_{N}\otimes\mathbf{a}^{T}\left(\theta\right)\right)\mathbf{x}\right)\\
= & \mathbf{x}^{H}\left(\mathbf{I}_{N}\otimes\mathbf{a}^{\ast}\left(\theta\right)\mathbf{a}^{T}\left(\theta\right)\right)\mathbf{x}=\mathbf{x}^{H}\mathbf{A}\left(\theta\right)\mathbf{x},
\end{aligned}
\]
where $\mathbf{A}\left(\theta\right)=\mathbf{I}_{N}\otimes\mathbf{a}^{\ast}\left(\theta\right)\mathbf{a}^{T}\left(\theta\right)$.

Suppose there are $\overline{K}$ targets of interest, and then the
spatial cross-correlation sidelobes (cross-correlation beampattern)
between the probing signals at locations $\theta_{i}$ and $\theta_{j}$
($i\neq j$, $i,j=1,\ldots,\overline{K}$ and $\theta_{i},\theta_{j}\in\Theta$)
is given by
\[
\begin{aligned} & P_{cc}\left(\theta_{i},\theta_{j},\mathbf{x}\right)\\
= & \sum_{n=1}^{N}\left(\mathbf{a}^{T}\left(\theta_{i}\right)\mathbf{x}\left(n\right)\right)^{\ast}\left(\mathbf{a}^{T}\left(\theta_{j}\right)\mathbf{x}\left(n\right)\right)\\
= & \left(\left(\mathbf{I}_{N}\otimes\mathbf{a}^{T}\left(\theta_{i}\right)\right)\mathbf{x}\right)^{H}\left(\left(\mathbf{I}_{N}\otimes\mathbf{a}^{T}\left(\theta_{j}\right)\right)\mathbf{x}\right)\\
= & \mathbf{x}^{H}\left(\mathbf{I}_{N}\otimes\mathbf{a}^{\ast}\left(\theta_{i}\right)\mathbf{a}^{T}\left(\theta_{j}\right)\right)\mathbf{x}=\mathbf{x}^{H}\mathbf{A}\left(\theta_{i},\theta_{j}\right)\mathbf{x},
\end{aligned}
\]
where $\mathbf{A}\left(\theta_{i},\theta_{j}\right)=\mathbf{I}_{N}\otimes\mathbf{a}^{\ast}\left(\theta_{i}\right)\mathbf{a}^{T}\left(\theta_{j}\right)$.

The objective of the transmit beampattern matching problem is as follows:
i) to match a desired transmit beampattern denoted as $p\left(\theta\right)$,
which can be formulated as follows\footnote{Variable $\alpha$ is introduced since $p\left(\theta\right)$ is
typically given in a \textquotedblleft normalized form\textquotedblright{}
and we want to approximate a scaled version of $p\left(\theta\right)$,
not $p\left(\theta\right)$ itself.}:
\begin{equation}
J\left(\alpha,\mathbf{x}\right)=\sum_{\theta\in\Theta}\omega\left(\theta\right)\left|\alpha p\left(\theta\right)-P\left(\theta,\mathbf{x}\right)\right|^{2},\label{eq:beampattern matching}
\end{equation}
where $\omega\left(\theta\right)\geq0$ is the weight for the direction
$\theta$; and ii) to minimize the cross-correlation between the probing
signals at a number of given target locations due to the fact that
the statistical performance of adaptive MIMO radar techniques rely
on the cross-correlation beampattern, which is given as
\begin{equation}
E\left(\mathbf{x}\right)=\sum_{\theta_{i},\theta_{j}\in\overline{\Theta},\,i\neq j}\left|P_{cc}\left(\theta_{i},\theta_{j},\mathbf{x}\right)\right|^{2}.\label{eq:sidelobe term}
\end{equation}
Then, by considering $J\left(\alpha,\mathbf{x}\right)$ and $E\left(\mathbf{x}\right)$,
the MIMO transmit beampattern matching problem is formulated as follows:
\begin{equation}
\begin{aligned} & \underset{\alpha,\mathbf{x}}{\mathsf{minimize}} &  & f\left(\alpha,\mathbf{x}\right)\triangleq J\left(\alpha,\mathbf{x}\right)+\omega_{cc}E\left(\mathbf{x}\right)\\
 & \mathsf{subject\:to} &  & \mathbf{x}\in\mathcal{X}\triangleq\mathcal{X}_{0}\cap\left(\cap_{i}\mathcal{X}_{i}\right),
\end{aligned}
\label{eq:problem}
\end{equation}
where $\omega_{cc}$ controls the sidelobe term, $\mathcal{X}$ generally
denotes the waveform constraint, and $\mathcal{X}_{0}=\left\{ \mathbf{x}\in\mathbb{C}^{MN}\mid\left\Vert \mathbf{x}\right\Vert _{2}^{2}=c_{e}^{2}\right\} $
representing the \textbf{total transmit energy (power) constraint}.
We are also interested in other practical waveform constraints:

\textbf{i) Constant modulus constraint }is to prevent the non-linearity
distortion of the power amplifier to maximize the efficiency of the
transmitter, which is given by $\begin{array}{c}
\mathcal{X}_{1}=\left\{ \mathbf{x}\mid\left|x\left(l\right)\right|=c_{d}=\frac{c_{e}}{\sqrt{MN}}\right\} \end{array}$ for $l=1,\ldots,MN$;

\textbf{ii) Peak-to-Average Ratio (PAR) constraint} is the ratio of
the peak signal power to its average power ($\mathrm{PAR}\left(\mathbf{x}\right)=\frac{\max\left|x\left(l\right)\right|^{2}}{\left\Vert \mathbf{x}\right\Vert _{2}^{2}/MN}$
with $1\leq\mathrm{PAR}\left(\mathbf{x}\right)\leq MN$). The $\mathrm{PAR}\left(\mathbf{x}\right)$
is constrained to a small threshold, so that the analog-to-digital
and digital-to-analog converters can have lower dynamic range, and
fewer linear power amplifiers are needed. Since $\mathcal{X}_{0}$,
the PAR constraint is $\mathcal{X}_{2}=\left\{ \mathbf{x}\mid\left|x\left(l\right)\right|\leq c_{p},\frac{c_{e}}{\sqrt{MN}}\leq c_{p}\leq c_{e}\right\} $
for $l=1,\ldots,MN$;

\textbf{iii) Similarity constraint} is to allow the designed waveforms
to lie in the neighborhood of a reference one which already can attain
a good performance \cite{LiGuerciXu2006a}, which is denoted as $\begin{array}{c}
\mathcal{X}_{3}=\left\{ \mathbf{x}\mid\left|\mathbf{x}-\mathbf{x}_{\mathrm{ref}}\right|\leq c_{\epsilon},0\leq c_{\epsilon}\leq\frac{2}{\sqrt{MN}}\right\} \end{array}$.

Problem \eqref{eq:problem} is a constrained nonconvex problem due
to the nonconvex objective and constraints. We are trying to solve
it by using efficient nonconvex optimization methods.

\section{Problem Solving via The MM Method}

\subsection{The Majorization-Minimization (MM) Method}

The MM method \cite{HunterLange2004,RazaviyaynHongLuo2013,SunBabuPalomar2016}
is a generalization of the well-known EM method. For an optimization
problem given by
\[
\begin{aligned} & \underset{\mathbf{x}}{\mathsf{minimize}} &  & f\left(\mathbf{x}\right)\\
 & \mathsf{subject\:to} &  & \mathbf{x}\in{\cal X},
\end{aligned}
\]
instead of dealing with this problem directly which could be difficult,
the MM-based algorithm solves a series of simpler subproblems with
surrogate functions that majorize $f\left(\mathbf{x}\right)$ over
${\cal X}$. More specifically, starting from an initial point $\mathbf{x}^{\left(0\right)}$,
it produces a sequence $\left\{ \mathbf{x}^{\left(k\right)}\right\} $
by the following update rule:
\[
\mathbf{x}^{\left(k\right)}\in\arg\min_{\mathbf{x}\in{\cal X}}\:\overline{f}\left(\mathbf{x},\mathbf{x}^{\left(k-1\right)}\right),
\]
where the surrogate majorizing function $\overline{f}\left(\mathbf{x},\mathbf{x}^{\left(k\right)}\right)$
satisfies
\[
\begin{array}{cl}
\overline{f}\left(\mathbf{x}^{\left(k\right)},\mathbf{x}^{\left(k\right)}\right)=f\left(\mathbf{x}^{\left(k\right)}\right), & \forall\mathbf{x}^{\left(k\right)}\in{\cal X},\\
\overline{f}\left(\mathbf{x},\mathbf{x}^{\left(k\right)}\right)\geq f\left(\mathbf{x}\right), & \forall\mathbf{x},\mathbf{x}^{\left(k\right)}\in{\cal X},\\
\overline{f}^{\prime}\left(\mathbf{x}^{\left(k\right)},\mathbf{x}^{\left(k\right)};\mathbf{d}\right)=f^{\prime}\left(\mathbf{x}^{\left(k\right)};\mathbf{d}\right), & \forall\mathbf{d},\mbox{\text{ s.t. }}\mathbf{x}^{\left(k\right)}+\mathbf{d}\in{\cal X}.
\end{array}
\]

The objective function value is monotonically nonincreasing at each
iteration. To use the MM method, the key step is to find a majorizing
function to make the subproblem easy to solve, which will be discussed
in the following subsections.

\subsection{Majorization Steps For The Beampattern Matching Term $J\left(\alpha,\mathbf{x}\right)$}

In this section, we discuss the majorization steps, i.e., how to construct
a good majorizing function for the beampattern matching term $J\left(\alpha,\mathbf{x}\right)$
in \eqref{eq:beampattern matching}. First, we have
\[
\begin{aligned}J\left(\alpha,\mathbf{x}\right)= & \sum_{\theta\in\Theta}\omega\left(\theta\right)\left|\alpha p\left(\theta\right)-P\left(\theta,\mathbf{x}\right)\right|^{2}\\
= & \alpha^{2}\sum_{\theta\in\Theta}\omega\left(\theta\right)p^{2}\left(\theta\right)-2\alpha\sum_{\theta\in\Theta}\omega\left(\theta\right)p\left(\theta\right)P\left(\theta,\mathbf{x}\right)+\sum_{\theta\in\Theta}\omega\left(\theta\right)\left(P\left(\theta,\mathbf{x}\right)\right)^{2},
\end{aligned}
\]
which is a quadratic function in variable $\alpha$. Then, it follows
that the minimum of $J\left(\alpha,\mathbf{x}\right)$ is attained
when
\[
\alpha\left(\mathbf{x}\right)=\sum_{\theta\in\Theta}\omega\left(\theta\right)p\left(\theta\right)P\left(\theta,\mathbf{x}\right)/\sum_{\theta\in\Theta}\omega\left(\theta\right)p^{2}\left(\theta\right).
\]
Substituting $\alpha\left(\mathbf{x}\right)$ back into $J\left(\alpha,\mathbf{x}\right)$
and considering
\[
P\left(\theta,\mathbf{x}\right)=\mathrm{Tr}\left(\mathbf{x}\mathbf{x}^{H}\mathbf{A}\left(\theta\right)\right)=\mathrm{vec}\left(\mathbf{x}\mathbf{x}^{H}\right)^{H}\mathrm{vec}\left(\mathbf{A}\left(\theta\right)\right),
\]
we get
\[
\begin{aligned}J\left(\mathbf{x}\right)= & \sum_{\theta\in\Theta}\omega\left(\theta\right)\left(\mathrm{vec}\left(\mathbf{x}\mathbf{x}^{H}\right)^{H}\mathrm{vec}\left(\mathbf{A}\left(\theta\right)\right)\right)^{2}-\left(\sum_{\theta\in\Theta}\omega\left(\theta\right)p^{2}\left(\theta\right)\right)^{-1}\\
 & \times\left(\mathrm{vec}\left(\mathbf{x}\mathbf{x}^{H}\right)^{H}\right.\mathrm{vec}\Bigl(\sum_{\theta\in\Theta}\omega\left(\theta\right)\left.p\left(\theta\right)\mathbf{A}\left(\theta\right)\Bigr)\right)^{2}\\
= & \mathrm{vec}\left(\mathbf{x}\mathbf{x}^{H}\right)^{H}\mathbf{H}_{J}\mathrm{vec}\left(\mathbf{x}\mathbf{x}^{H}\right),
\end{aligned}
\]
where
\[
\begin{aligned}\mathbf{H}_{J}= & \sum_{\theta\in\Theta}\omega\left(\theta\right)\mathrm{vec}\left(\mathbf{A}\left(\theta\right)\right)\mathrm{vec}\left(\mathbf{A}\left(\theta\right)\right)^{H}-\Bigl(\sum_{\theta\in\Theta}\omega\left(\theta\right)p^{2}\left(\theta\right)\Bigr)^{-1}\\
 & \times\mathrm{vec}\Bigl(\sum_{\theta\in\Theta}\omega\left(\theta\right)p\left(\theta\right)\mathbf{A}\left(\theta\right)\Bigr)\mathrm{vec}\Bigl(\sum_{\theta\in\Theta}\omega\left(\theta\right)p\left(\theta\right)\mathbf{A}\left(\theta\right)\Bigr)^{H},
\end{aligned}
\]
and it is easy to see that $J\left(\mathbf{x}\right)$ is a quartic
function in $\mathbf{x}$. Next, we introduce a useful lemma.
\begin{lem}
\label{lem:quadratic majorization} Let $\mathbf{A}\in\mathbb{H}^{K}$
and $\mathbf{B}\in\mathbb{H}^{K}$ such that $\mathbf{B}\succeq\mathbf{A}$.
At any point $\mathbf{x}_{0}\in\mathbb{C}^{K}$, the quadratic function
\textbf{$\mathbf{x}^{T}\mathbf{A}\mathbf{x}$} is majorized by $\mathbf{x}^{H}\mathbf{B}\mathbf{x}+2\mathrm{Re}\left(\mathbf{x}^{H}\left(\mathbf{A}-\mathbf{B}\right)\mathbf{x}_{0}\right)+\mathbf{x}_{0}^{H}\left(\mathbf{B}-\mathbf{A}\right)\mathbf{x}_{0}$.
\end{lem}
\begin{IEEEproof}
Notice that $\left(\mathbf{x}-\mathbf{x}_{0}\right)^{H}\left(\mathbf{B}-\mathbf{A}\right)\left(\mathbf{x}-\mathbf{x}_{0}\right)\geq0$.

Based on Lemma \ref{lem:quadratic majorization}, we can choose $\psi_{J,1}\geq\lambda_{\mathrm{max}}\left(\mathbf{H}_{J}\right)$,
and because $\psi_{J,1}\mathbf{I}\succeq\mathbf{H}_{J}$, at iterate
$\mathbf{x}^{\left(t\right)}$ we have
\[
\begin{aligned}J\left(\mathbf{x}\right)\leq & \psi_{J,1}\mathrm{vec}\left(\mathbf{x}\mathbf{x}^{H}\right)^{H}\mathrm{vec}\left(\mathbf{x}\mathbf{x}^{H}\right)\\
 & +2\mathrm{Re}\left(\mathrm{vec}\left(\mathbf{x}\mathbf{x}^{H}\right)^{H}\left(\mathbf{H}_{J}-\psi_{J,1}\mathbf{I}\right)\mathrm{vec}\left(\mathbf{x}^{\left(t\right)}\mathbf{x}^{\left(t\right)H}\right)\right)\\
 & +\mathrm{vec}\left(\mathbf{x}^{\left(t\right)}\mathbf{x}^{\left(t\right)H}\right)^{H}\left(\psi_{J,1}\mathbf{I}-\mathbf{H}_{J}\right)\mathrm{vec}\left(\mathbf{x}^{\left(t\right)}\mathbf{x}^{\left(t\right)H}\right),
\end{aligned}
\]
where since $\mathrm{vec}\left(\mathbf{x}\mathbf{x}^{H}\right)^{H}\mathrm{vec}\left(\mathbf{x}\mathbf{x}^{H}\right)=\left\Vert \mathbf{x}\right\Vert _{2}^{4}=c_{e}^{4}$,
the first term is just a constant. Then after ignoring the constant
terms, we get the following majorizing function for $J\left(\mathbf{x}\right)$:
\[
\overline{J}_{1}\left(\mathbf{x},\mathbf{x}^{\left(t\right)}\right)\simeq2\mathrm{Re}\left(\mathrm{vec}\left(\mathbf{x}\mathbf{x}^{H}\right)^{H}\left(\mathbf{H}_{J}-\psi_{J,1}\mathbf{I}\right)\mathrm{vec}\left(\mathbf{x}^{\left(t\right)}\mathbf{x}^{\left(t\right)H}\right)\right),
\]
where ``$\simeq$'' stands for ``equivalence'' up to additive
constants. Substituting $\mathbf{H}_{J}$ back into function $\overline{J}_{1}\left(\mathbf{x},\mathbf{x}^{\left(t\right)}\right)$
and dropping the constants, we have
\begin{equation}
\overline{J}_{1}\left(\mathbf{x},\mathbf{x}^{\left(t\right)}\right)\simeq2\mathbf{x}^{H}\left(\mathbf{M}_{J}-\psi_{J,1}\mathbf{x}^{\left(t\right)}\mathbf{x}^{\left(t\right)H}\right)\mathbf{x},
\end{equation}
where $\mathbf{M}_{J}=$$\sum_{\theta\in\Theta}\omega\left(\theta\right)\left(P\left(\theta,\mathbf{x}^{\left(t\right)}\right)-p\left(\theta\right)\alpha\left(\mathbf{x}^{\left(t\right)}\right)\right)\mathbf{A}\left(\theta\right)$.
It is easy to see that after majorization, the majorizing function
$\overline{J}_{1}\left(\mathbf{x},\mathbf{x}^{\left(t\right)}\right)$
becomes quadratic in $\mathbf{x}$ rather than quartic in $J\left(\mathbf{x}\right)$.
However, using this function as the objective to solve is still hard
due to the waveform constraint ${\cal X}$.\footnote{It is a NP-hard unimodular quadratic program even only considering
$\mathcal{X}_{1}$.} So we propose to majorize $\overline{J}_{1}\left(\mathbf{x},\mathbf{x}^{\left(t\right)}\right)$
again to simplify the problem to solve in each iteration. Thus, we
can consider choosing $\psi_{J,2}\geq\lambda_{\mathrm{max}}\left(\mathbf{M}_{J}\right)\geq\lambda_{\mathrm{max}}\left(\mathbf{M}_{J}-\psi_{J,1}\mathbf{x}^{\left(t\right)}\mathbf{x}^{\left(t\right)H}\right)$
for majorization, where we can have the following useful property.
\end{IEEEproof}
\begin{lem}
\label{lem:Hermitian Toeplitz}\cite{JorgeFerreira1994,HornJohnson1990}
Define 
\[
\begin{aligned}\mathbf{B}= & \sum_{\theta\in\Theta}\omega\left(\theta\right)\Bigl(P\left(\theta,\mathbf{x}^{\left(t\right)}\right)-p\left(\theta\right)\alpha\left(\mathbf{x}^{\left(t\right)}\right)\Bigr)\mathbf{a}^{\ast}\left(\theta\right)\mathbf{a}^{T}\left(\theta\right)\\
= & \left[\begin{array}{cccc}
b_{0} & b_{1}^{\ast} & \cdots & b_{M-1}^{\ast}\\
b_{1} & b_{0} & \ddots & \vdots\\
\vdots & \ddots & \ddots & b_{1}^{\ast}\\
b_{M-1} & \ldots & b_{1} & b_{0}
\end{array}\right],
\end{aligned}
\]
which is Hermitian Toeplitz, $\mathbf{F}$ as a $2M\times2M$ FFT
matrix, and $\mathbf{b}=\left[b_{0},b_{1},\ldots,b_{M-1},0,b_{M-1}^{\ast},,\ldots,b_{1}^{\ast}\right]^{T}$.
Then, we have $\mathbf{M}_{J}=\mathbf{I}_{N}\otimes\mathbf{B}$, $\lambda_{\mathrm{max}}\left(\mathbf{M}_{J}\right)=\lambda_{\mathrm{max}}\left(\mathbf{B}\right)$,
and
\[
\lambda_{\mathrm{max}}\left(\mathbf{B}\right)\leq\lambda_{\mu}=\frac{1}{2}\left(\underset{1\leq i\leq M}{\max}\mu_{2i}+\underset{1\leq i\leq M}{\max}\mu_{2i-1}\right),
\]
where $\boldsymbol{\mu}=\mathbf{F}\mathbf{b}$, which is the discrete
Fourier transform for $\mathbf{b}$.
\end{lem}
Lemma \ref{lem:Hermitian Toeplitz} provides an easy way for the computation
of $\psi_{J,2}.$ Based on Lemma \ref{lem:quadratic majorization}
and using $\psi_{J,2}=\lambda_{\mu}$, the majorizing function $\overline{J}_{1}\left(\mathbf{x},\mathbf{x}^{\left(t\right)}\right)$
can be further majorized as
\[
\begin{aligned}\overline{J}_{1}\left(\mathbf{x},\mathbf{x}^{\left(t\right)}\right)\leq & 2\psi_{J,2}\mathbf{x}^{H}\mathbf{x}+4\mathrm{Re}\Bigl(\mathbf{x}^{H}\Bigl(\mathbf{M}_{J}-\psi_{J,1}\mathbf{x}^{\left(t\right)}\mathbf{x}^{\left(t\right)H}-\psi_{J,2}\mathbf{I}\Bigr)\mathbf{x}^{\left(t\right)}\Bigr)\\
 & +2\mathbf{x}^{\left(t\right)H}\left(\psi_{J,2}\mathbf{I}-\mathbf{M}_{J}+\psi_{J,1}\mathbf{x}^{\left(t\right)}\mathbf{x}^{\left(t\right)H}\right)\mathbf{x}^{\left(t\right)},
\end{aligned}
\]
where since $\left\Vert \mathbf{x}\right\Vert _{2}^{2}=c_{e}^{2}$,
the first term is a constant. Then by ignoring the constant terms,
the objective becomes a linear majorizing function at iterate $\mathbf{x}^{\left(t\right)}$
as follows:
\begin{equation}
\overline{J}_{2}\left(\mathbf{x},\mathbf{x}^{\left(t\right)}\right)\simeq-4\mathrm{Re}\left(\mathbf{x}^{H}\mathbf{y}_{J}\right),
\end{equation}
where $\mathbf{y}_{J}=-\left(\mathbf{M}_{J}-c_{e}^{2}\psi_{J,1}\mathbf{I}-\psi_{J,2}\mathbf{I}\right)\mathbf{x}^{\left(t\right)}$.

\subsection{Majorization Steps For The Sidelobe Term $E\left(\mathbf{x}\right)$}

To deal with the sidelobe term $E\left(\mathbf{x}\right)$ in \eqref{eq:sidelobe term},
the majorization steps are similar to $J\left(\mathbf{x}\right)$.
First, we have
\[
\begin{aligned}E\left(\mathbf{x}\right)= & \sum_{\theta_{i},\theta_{j}\in\overline{\Theta},\,i\neq j}\left|P_{cc}\left(\theta_{i},\theta_{j},\mathbf{x}\right)\right|^{2}\\
= & \mathrm{vec}\left(\mathbf{x}\mathbf{x}^{H}\right)^{H}\mathbf{H}_{E}\mathrm{vec}\left(\mathbf{x}\mathbf{x}^{H}\right),
\end{aligned}
\]
where $\mathbf{H}_{E}=\sum_{\theta_{i},\theta_{j}\in\overline{\Theta},\,i\neq j}\mathrm{vec}\left(\mathbf{A}\left(\theta_{i},\theta_{j}\right)\right)\mathrm{vec}\left(\mathbf{A}\left(\theta_{i},\theta_{j}\right)\right)^{H}$.
Then, based on Lemma \ref{lem:quadratic majorization}, by choosing
$\psi_{E,1}\geq\lambda_{\mathrm{max}}\left(\mathbf{H}_{E}\right)$
and $\psi_{E,2}\geq\lambda_{\mathrm{max}}\left(\mathbf{M}_{E}-\psi_{E,1}\mathbf{x}^{\left(t\right)}\mathbf{x}^{\left(t\right)H}\right)$,
we can get the majorizing functions at iterate $\mathbf{x}^{\left(t\right)}$
written as follows:
\begin{equation}
\begin{aligned}\overline{E}_{1}\left(\mathbf{x},\mathbf{x}^{\left(t\right)}\right)\simeq & 2\mathbf{x}^{H}\left(\mathbf{M}_{E}-\psi_{E,1}\mathbf{x}^{\left(t\right)}\mathbf{x}^{\left(t\right)H}\right)\mathbf{x}\\
\leq & \overline{E}_{2}\left(\mathbf{x},\mathbf{x}^{\left(t\right)}\right)\\
\simeq & -4\mathrm{Re}\left(\mathbf{x}^{H}\mathbf{y}_{E}\right),
\end{aligned}
\end{equation}
where $\mathbf{M}_{E}=\sum_{\theta_{i},\theta_{j}\in\overline{\Theta},\,i\neq j}P_{cc}\left(\theta_{j},\theta_{i},\mathbf{x}^{\left(t\right)}\right)\mathbf{A}\left(\theta_{i},\theta_{j}\right)$
and $\mathbf{y}_{E}=-\left(\mathbf{M}_{E}-c_{e}^{2}\psi_{E,1}\mathbf{I}-\psi_{E,2}\mathbf{I}\right)\mathbf{x}^{\left(t\right)}$.

\subsection{Solving The Majorized Subproblem in MM}

By combing the two majorizing functions $\overline{J}_{2}\left(\mathbf{x},\mathbf{x}^{\left(t\right)}\right)$
and $\overline{E}_{2}\left(\mathbf{x},\mathbf{x}^{\left(t\right)}\right)$,
the overall majorizing function at iterate $\mathbf{x}^{\left(k\right)}$
for the objective $f\left(\mathbf{x}\right)$ is given as follows:
\[
\begin{aligned}f\left(\mathbf{x}\right)\leq & \overline{f}\left(\mathbf{x},\mathbf{x}^{\left(t\right)}\right)\\
= & \overline{J}_{2}\left(\mathbf{x},\mathbf{x}^{\left(t\right)}\right)+\omega_{cc}\overline{E}_{2}\left(\mathbf{x},\mathbf{x}^{\left(t\right)}\right)\\
\simeq & -4\mathrm{Re}\left(\mathbf{x}^{H}\mathbf{y}_{J}\right)-4\omega_{cc}\mathrm{Re}\left(\mathbf{x}^{H}\mathbf{y}_{E}\right)\\
= & -\mathrm{Re}\left(\mathbf{x}^{H}\mathbf{y}\right),
\end{aligned}
\]
where
\[
\begin{aligned}\mathbf{y}= & -4\left(\mathbf{M}_{J}+\omega_{cc}\mathbf{M}_{E}-c_{e}^{2}\left(\psi_{J,1}+\omega_{cc}\psi_{E,1}\right)\mathbf{I}\right.\\
 & -\left.\left.\left(\psi_{J,2}\right.+\omega_{cc}\psi_{E,2}\right)\mathbf{I}\right)\mathbf{x}^{\left(t\right)}.
\end{aligned}
\]

Finally, by majorizing the objective function in \eqref{eq:problem}
using the MM method, the subproblem we need to solve at each iteration
is given as follows:
\begin{equation}
\begin{aligned} & \mathrm{minimize}_{\mathbf{x}} &  & \overline{f}\left(\mathbf{x},\mathbf{x}^{\left(t\right)}\right)\simeq-\mathrm{Re}\left(\mathbf{x}^{H}\mathbf{y}\right)\\
 & \mathrm{subject\:to} &  & \mathbf{x}\in{\cal X}.
\end{aligned}
\label{eq:subproblem}
\end{equation}
For problem \eqref{eq:subproblem}, as to different interested waveform
constraints, closed-form optimal solutions $\mathbf{x}^{\star}$ can
be derived, which are summarized in the following lemma.
\begin{lem}
\label{lem:closed-form solution}\textbf{i)} For fixed energy constraint
(i.e., ${\cal X}={\cal X}_{0}$), $\mathbf{x}^{\star}=c_{e}\mathbf{y}/\left\Vert \mathbf{y}\right\Vert _{2}$;
\textbf{ii)} for constant modulus constraint (i.e., ${\cal X}={\cal X}_{1}$),
$\mathbf{x}^{\star}=c_{d}e^{j\arg\left(\mathbf{y}\right)}$;\footnote{The operation $\arg\left(\mathbf{y}\right)$ is applied element-wise
for $\mathbf{y}$.} \textbf{iii)} for fixed energy with PAR constraint (i.e., ${\cal X}={\cal X}_{0}\cap{\cal X}_{2}$),
the solution $\mathbf{x}^{\star}$ can be found in \cite[Alg. 2]{TroppDhillonHeathStrohmer2005};
\textbf{iv)} for constant modulus with similarity constraint (i.e.,
${\cal X}={\cal X}_{1}\cap{\cal X}_{3}$), the solution $\mathbf{x}^{\star}$
can be found in \cite{ZhaoPalomar2017}.
\end{lem}

\subsection{The MM-Based Beampattern Matching Algorithm}

Based on the MM method, in order to solve the original problem \eqref{eq:problem},
we just need to iteratively solve the subproblem \eqref{eq:subproblem}
with a closed-form solution update in Lemma \ref{lem:closed-form solution}
at each iteration. The overall algorithm is summarized as follows.

\noindent\fbox{\begin{minipage}[t]{1\columnwidth - 2\fboxsep - 2\fboxrule}%
\textbf{Input:} $\mathbf{a}\left(\theta\right)$, $p\left(\theta\right)$,\textbf{
}$\mathbf{x}^{\left(0\right)}$ and $t=0$. 

\textbf{Repeat}

$\:\:$1. Compute $\mathbf{M}_{J}$, $\mathbf{M}_{E},$ $\psi_{J,1}$,
$\psi_{E,1}$, $\psi_{J,2}$, $\psi_{E,2}$ and $\mathbf{y}$;

$\:\:$2. Update $\mathbf{x}^{\left(t\right)}$ in a closed-form according
to Lemma \ref{lem:closed-form solution};

$\:\:$3. $t=t+1$;

\textbf{Until} $\mathbf{x}$ and $f\left(\mathbf{x}\right)$ satisfy
a termination criterion.

\textbf{Output:} $\alpha$, $\mathbf{x}$.%
\end{minipage}}

\section{Numerical Simulations\label{sec:Numerical-Simulations}}

The performance of the proposed algorithm for MIMO transmit beampattern
matching is evaluated by numerical simulations. A colocated MIMO radar
system is considered with a ULA comprising $M=10$ antennas with half-wavelength
spacing between adjacent antennas. Without loss of generality, the
total transmit power is set to $c_{e}^{2}=1$. Each transmit pulse
has $N=32$ samples. The range of angle is $\Theta=\left(-90^{\circ},90^{\circ}\right)$
with spacing $1^{\circ}$ under which the weight $\omega\left(\theta\right)=1$
for $\theta\in\Theta$, and $\omega_{cc}=0$, which is the same setting
as \cite{ChengHeZhangLi2017}. We consider a desired beampattern with
three targets or mainlobes ($K=3$) at $\theta_{1}=-40^{\circ}$,
$\theta_{2}=0^{\circ}$, $\theta_{3}=40^{\circ}$, and each width
of them is $\triangle\theta=20^{\circ}$. The desired beampattern
is
\[
p\left(\theta\right)=\begin{cases}
1, & \theta\in\left[\theta_{k}-\triangle\theta/2,\theta_{k}+\triangle\theta/2\right],\,k=1,2,K\\
0, & \mathrm{otherwise}.
\end{cases}
\]

\begin{figure}[H]
\begin{centering}
\includegraphics[scale=0.43]{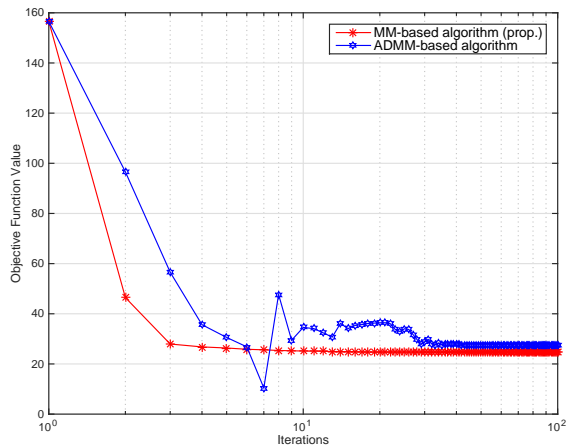}
\par\end{centering}
\centering{}\caption{\label{fig:Convergence-comparison-for}Convergence comparison for
objective function value.}
\end{figure}

We compare the convergence property over iterations of the objective
function for the beampattern matching problem under unimodulus waveform
constraint by using the proposed MM-based algorithm (denoted as MM-based
algorithm (prop.)) and the ADMM-based algorithm in \cite{ChengHeZhangLi2017}
(denoted as ADMM-based algorithm) , which is shown in Fig. \ref{fig:Convergence-comparison-for}.

As shown in Fig. \ref{fig:Convergence-comparison-for}, the MM-based
algorithm can have a monotonic convergence property. And it can converge
within $20$ iterations which is faster than the benchmark algorithm. 

\begin{figure}[H]
\begin{centering}
\includegraphics[scale=0.43]{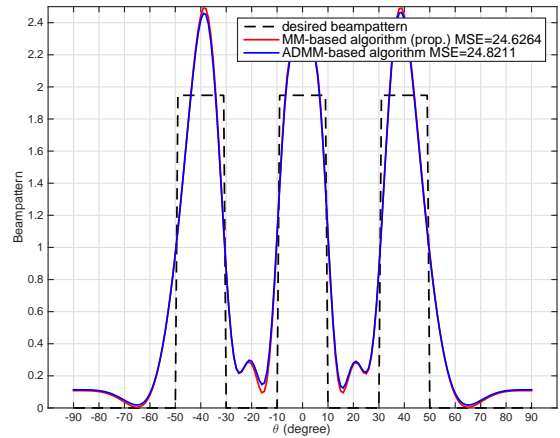}
\par\end{centering}
\caption{\label{fig:beampattern-matching}Transmit beampattern design with
$3$ targets .}
\end{figure}

Then, we also compare the matching performance of the designed beampatterns
in terms of the mean-squared error (MSE) defined as
\[
\text{MSE}\left(P\left(\theta,\mathbf{x}\right)\right)=\mathbb{E}\left[\sum_{\theta\in\Theta}\omega\left(\theta\right)\left|\alpha p\left(\theta\right)-P\left(\theta,\mathbf{x}\right)\right|^{2}\right].
\]
In Fig. \ref{fig:beampattern-matching}, we show the simulation results
for $\text{MSE}\left(P\left(\theta,\mathbf{x}\right)\right)$ by using
different design methods.

From Fig. \ref{fig:beampattern-matching}, we can see that compared
to the benchmark, our proposed algorithm can have a tighter matching
performance and can obtain a lower MSE. Based on these, the proposed
algorithm is validated.

\section{Conclusions\label{sec:Conclusions}}

This paper has considered the MIMO transmit beampattern matching problem.
Efficient algorithms have been proposed based on the MM method. Numerical
simulations show that the proposed algorithms are efficient in solving
the beampattern matching problem and can obtain a better performance
compared to the the state-of-art method.

\bibliographystyle{IEEEtran}
\bibliography{/Users/ziping/Dropbox/Research/1-Report/Reference/radar}

\end{document}